\documentclass[11pt]{article}
\usepackage{mathrsfs}
\usepackage{amsfonts}
\usepackage{amssymb}
\usepackage{graphicx}
\usepackage{amsmath}

\def\sqr#1#2{{\vcenter{\vbox{\hrule height.#2pt
              \hbox{\vrule width.#2pt height#1pt \kern#1pt \vrule
width.#2pt}
              \hrule height.#2pt}}}}
\def\signed #1{{\unskip\nobreak\hfil\penalty50
              \hskip2em\hbox{}\nobreak\hfil#1
              \parfillskip=0pt \finalhyphendemerits=0 \par}}
\def\endpf{\signed {$\sqr69$}}
\def\3n{\negthinspace \negthinspace \negthinspace }
\def\2n{\negthinspace \negthinspace }
\def\1n{\negthinspace }

\def\ms{\medskip}

\def\({\Big (}
\def\){\Big )}
\def\[{\Big[}
\def\]{\Big]}

\def\be{\begin{equation}}
\def\bel{\begin{equation}\label}
\def\ee{\end{equation}}
\def\bea{\begin{eqnarray}}
\def\eea{\end{eqnarray}}
\def\bt{\begin{theorem}}
\def\et{\end{theorem}}
\def\bc{\begin{corollary}}
\def\ec{\end{corollary}}
\def\bl{\begin{lemma}}
\def\el{\end{lemma}}
\def\bp{\begin{proposition}}
\def\ep{\end{proposition}}
\def\br{\begin{remark}}
\def\er{\end{remark}}
\def\ba{\begin{array}}
\def\ea{\end{array}}
\def\bd{\begin{definition}}
\def\ed{\end{definition}}

\newtheorem{lemma}{Lemma}[section]
\newtheorem{remark}{Remark}[section]

\newtheorem{theorem}{Theorem}[section]
\newtheorem{corollary}{Corollary}[section]

\newtheorem{definition}{Definition}[section]
\newtheorem{proposition}{Proposition}[section]

\makeatletter
   
   \@addtoreset{equation}{section}
\makeatother
\allowdisplaybreaks
\begin{document}

\title{\bf Logarithmic wave equations in non-cylindrical domains\thanks{This work is
supported by the NSF of China under grants 11471070, 11771074 and
11371084.}}

\author{Lingyang Liu\thanks{School of Mathematics and Statistics, Northeast Normal
University, Changchun 130024, China. E-mail address:
liuly938@nenu.edu.cn. \ms }
}

\date{}

\maketitle

\begin{abstract}
This paper is devoted to studying a type of logarithmic wave equation in non-cylindrical domains. Firstly, by the penalty method, we prove the existence of weak solutions to such kind of equations. Secondly, different from the dissipative wave equation, the energy defined in this problem is not always positive. Thus, some suitable initial data are selected to let the energy be positive. Finally, by a difference inequality, we derive a exponential decay estimate for the positive energy.

\end{abstract}

\noindent{\bf Key Words. } logarithmic wave equation, weak solution, energy decay, non-cylindrical domain.

\section{Introduction and main results}

Let $T>0$, and $\Omega_t$ be a bounded domain in $\mathbb R,$ for $t\in[0,T].$ Set $\widehat{Q}_T=\Omega_t\times(0,T),$ and denote by $\widehat{\Sigma}_T$ the lateral boundary of $\widehat{Q}_T.$
Consider the following semi-linear wave equation with a logarithmic nonlinearity in the non-cylindrical domain $\widehat{Q}_T$:
\begin{equation}\label{1.1}
\left\{\begin{array}{ll}
u^{''}-\Delta u+au'+bu=u\ln|u|^\gamma& (x,t)\in \widehat{Q}_T,\\[2mm]
u=0 &(x,t)\in\widehat{\Sigma}_T,\\[2mm]
u(x,0)=u_0(x), \ u'(x,0)=u_1(x)&x\in\Omega_0,
\end{array}\right.
\end{equation}
where $(u_0, u_1)$ is any given initial couple, $(u, u')$ is the state variable, $a,b>0,$ and $0<\gamma<1.$

First, we give the definition of weak solutions to (\ref{1.1}):

\begin{definition}\label{d1.2} a function $u$ is called a weak solution of (\ref{1.1}), if
 $u\in L^\infty\big(0,T;H^1_0(\Omega_t)\big),$ $u'\in L^\infty\big(0,T;L^2(\Omega_t)\big),$
and it satisfies the integral equality
\begin{equation*}
-\int_{\widehat{Q}_T}u'\phi'dxdt+\int_{\widehat{Q}_T}\nabla u \nabla\phi dxdt+\int_{\widehat{Q}_T}\big(au'+bu\big)\phi dxdt=\int_{\widehat{Q}_T} u\ln|u|^\gamma\phi dxdt,
\end{equation*}
for all $\phi$ such that $\phi\in L^1\big(0,T; H^1_0(\Omega_t)\big),$ $\phi'\in L^1\big(0,T; L^2(\Omega_t)\big)$ with $\phi(x,0)=\phi(x,T)=0.$
\end{definition}

There have been numerous works concerned with nonlinear hyperbolic equations in cylindrical domains (see e.g.\cite{gor,han,hu,lia,LMF,zua} and rich references therein). Being one of the most important nonlinear hyperbolic equations, logarithmic wave equations have received more and more attention. For example, \cite{lia} considered a class of logarithmic wave equation and related well-posedness problem. \cite{hu} investigated the asymptotic behavior of a type of logarithmic wave equation with linear damping. We also refer to \cite{gor,han} for some known results in this respect. This type of nonlinearities appear naturally in inflation cosmology and supersymmetric field theory, and has numerous applications in many branches of physics such as nuclear physics, optics and geophysics. Although it has been studied extensively in cylindrical cases, to the best of our knowledge, there has no paper concerned with this kind of issue in non-cylindrical domains.

This paper is devoted to a study of logarithmic wave equations in non-cylindrical domains. On solvability of initial-boundary value problem for nonlinear hyperbolic equations in non-cylindrical domains, in general, there are two ways to be chosen. On the one hand, exploiting hyperbolic diffeomorphisms, one can convert such problems into equivalent cylindrical problems and then apply the Faedo-Galerkin method to obtain a regular solution. For example, in \cite{fer}, the following nonlinear hyperbolic-parabolic equation with Dirichlet boundary condition was considered.
\begin{equation*}
  K_1(x,t)u_{tt}+K_2(x,t)u_t-\Delta u+f_1(t)|u|^\rho u=f(x,t),\quad \mbox{in}\ \widehat{Q},
\end{equation*}
where $K_1,K_2,f_1$ are given functions and $\widehat{Q}$ is some non-cylindrical domain. Under some hypotheses, the global existence and uniqueness of regular solutions were proved in this work. In \cite{lar}, a nonlinear Neumann type condition was put on the moving boundary and the authors studied a general nonlinearity $\Phi(u),$ which satisfies
\begin{equation*}
  \Phi\in C^1(\mathbb R),\quad \mbox{and} \quad F(u)=\int^u_0\Phi(s)ds\geq0.
\end{equation*}
These have played active roles in the second-order estimate for solutions. In \cite{ha}, a damped Klein-Gordon equation in non-cylindrical domains was investigated. By controlling the geometry of non-cylindrical domains, a hypothesis that appears in \cite{fer} is removed.

We consider the logarithmic nonlinearity. It doesn't meet any of the conditions mentioned above. At the beginning, we want to prove the existence of regular solutions for such a problem. However, the hardship rising in this situation is to get an upper bounded estimate for the term: $\big((u\ln |u|^\gamma)', u''\big),$ which appears in the second-order estimate for the solution $u.$

As we know, the method of diffeomorphisms turns a simple wave equation (in non-cylindrical domains) into a complicated wave equation (with more complex variable coefficients in cylindrical domains), that takes us a large amount of calculations. Besides, in order to obtain a desired decay estimate, restrictions for the boundary conditions are sometimes necessary. Instead of exploiting the method of diffeomorphisms, we adopt a penalty method, which is first introduced by J.L.Lions in \cite{lio}, to look for a weak solution for (\ref{1.1}). The main idea on this technology is to continue the original non-cylindrical equation into an approximate cylindrical one and the key point of this method is to obtain an uniform estimate for approximate solutions (see e.g.\cite{coo,lio}). Later, this technique is also used to study the decay of solutions and it is sufficient to make an estimate for approximate solutions independent of parameters (see e.g.\cite{nak,rab}). The advantage of this method is that there are not many requirements on geometry of domains and the fundamental logarithmic Sobolev inequality developed in fixed area can be applied directly. Based on this method, we prove the existence of weak solutions for (\ref{1.1}).

\medskip

We define the energy of (\ref{1.1})
\begin{equation}\label{en}
\displaystyle E(t)=\frac{1}{2}\|u'(t)\|^2+\frac{1}{2}\|\nabla u(t)\|^2+\frac{b}{2}\|u(t)\|^2-\frac{1}{2}\int_\Omega u^2(t)\ln|u(t)|^\gamma dx+\frac{\gamma}{4}\|u(t)\|^2.
\end{equation}

\medskip

Unlike nonlinear dissipation, energy (\ref{en}) is not always positive for any given initial value. Hence, we need to define a set where the initial data is selected so that the energy is positive. Further, by a difference inequality, a decay estimate of solutions was derived.

\medskip

In order to study the well-posedness of (\ref{1.1}), we display the following hypotheses on domain $\widehat{Q}_T:$

\medskip

\noindent (H1) (Geometric hypothesis) The family $\{\Omega_t\}_{t\in[0,T]}$ is increasing, i.e., if $t_1\leq t_2,$  then $\Omega_{t_1}\subseteq\Omega_{t_2};$

\medskip

\noindent (H2) (Regularity hypothesis) If $v\in H_0^1(\Omega)$ for $\Omega\supseteq\Omega_t,$ and $v=0\ a.e. $ in $\Omega-\Omega_t,$ then $v\in H^1_0(\Omega_t).$

\medskip

Next, we introduce some useful concepts and lemmas.

\medskip

For any $u\in H^1_0(\Omega),$ set
\begin{equation*}
  I_1(u)=\|\nabla u\|^2-\int_\Omega u^2\ln|u|^\gamma dx,
\end{equation*}
\begin{equation*}
  J_1(u)=\frac{1}{2}\|\nabla u\|^2-\frac{1}{2}\int_\Omega u^2\ln|u|^\gamma dx+\frac{\gamma}{4}\int_\Omega u^2dx,
\end{equation*}
\begin{equation*}
  N=\{u\in H^1_0(\Omega)\mid I_1(u)=0,\ \int_\Omega|\nabla u|^2dx\neq0\},
\end{equation*}
\begin{equation*}
  d=\inf\{\sup\limits_{\lambda\geq0} J_1(\lambda u)\mid u\in H^1_0(\Omega),\ \int_\Omega|\nabla u|^2dx\neq0\}.
\end{equation*}

\begin{lemma}[\cite{CLL}]
For any $u\in H^1_0(\Omega),$ $\displaystyle\int_\Omega|\nabla u|^2dx\neq0,$ let $g(\lambda)=J_1(\lambda u),$ then
\begin{eqnarray*}
\begin{array}{rl}
I_1(\lambda u)=\lambda g'(\lambda)=\left\{
                                   \begin{array}{ll}
                                     >0, & 0<\lambda<\lambda^*, \\
                                     =0, & \lambda=\lambda^*,\\
                                     <0, &  \lambda^*<\lambda<+\infty,
                                   \end{array}
                                 \right.
\end{array}
\end{eqnarray*}
where
\begin{equation*}
 \displaystyle\lambda^*=\exp\Big(\frac{\|\nabla u\|^2-\int_\Omega u^2\ln|u|^\gamma dx}{\|u\|}\Big).
\end{equation*}
Thus, it is easy to see that
\begin{equation*}
 d=\inf\limits_{u\in N}J_1(u).
\end{equation*}
\end{lemma}

Furthermore, $d$ has a lower bound.

\begin{lemma}[\cite{lia}]
$\displaystyle d\geq\frac{1}{4}e(\frac{2\pi}{\gamma})^{\frac{1}{2}}.$
\end{lemma}

\medskip

The well-posedness result for the system (\ref{1.1}) is stated as follows.
\begin{theorem}\label{t1.3}
Given $T>0.$ Suppose that conditions (H1) and (H2) hold. Then for any $u_0\in H^1_0(\Omega_0),$ $ u_1\in L^2(\Omega_0),$
the system $(\ref{1.1})$ admits a solution
$u\in L^\infty\big(0,T;H^1_0(\Omega_t)\big)\cap W^{1,\infty}\big(0,T;L^2(\Omega_t)\big),$
in the sense of Definition \ref{d1.2}.
\end{theorem}

\medskip

The decay result of solutions to (\ref{1.1}) can be stated as follows.

\begin{theorem}\label{t1}
Further assume that $0<E(0)<d,$ and $I_1(u)>0.$ Then the solution $u$ in Theorem \ref{t1.3} satisfies the following decay property:
\begin{equation*}
  E(t)\leq E(0)e^{-\beta t},
\end{equation*}
where $\beta$ is a positive constant.
\end{theorem}

\begin{remark}
The conditions on decay of solutions proposed in this paper are more flexible than that in \cite{hu}, which requires $0<E(0)<c_\gamma<d,$ for some inserted constant $c_\gamma$ (related to parameter $\gamma$).
\end{remark}

The rest of this paper is organized as follows. In Section 2,
some notations and fundamental lemmas are presented. In Section 3, we are devoted to proving the existence of weak solutions for (\ref{1.1}) (Theorem {\ref{t1.3}}). Finally, in Section 4, a decay estimate for the energy is given (Theorem {\ref{t1}}).

\section{Notations and primary lemmas}

In this section, we present some required definitions and lemmas that will be applied in the next sections. First, we show the logarithmic Sobolev inequality and the logarithmic Gronwall inequality which are important tools to obtain some desired estimates in proving the wellposedness of solutions.

\begin{lemma}[\cite{lia}]\label{log}(Logarithmic Sobolev Inequality)
Let $\Omega\subset\mathbb R^n$ be a bounded domain. u is a function in $H^1_0(\Omega)$ and $a\in\mathbb R^+.$ Then
\begin{equation*}
  2\int_\Omega|u(x)|^2\ln\big(\frac{|u(x)|}{\|u(x)\|}\big)dx+n(1+\ln a)\|u(x)\|^2\leq \frac{a^2}{\pi}\int_\Omega|\nabla u|^2dx.
\end{equation*}
\end{lemma}

\begin{lemma}[\cite{han}]\label{log}(Logarithmic Gronwall Inequality)
Assume that w(t) is nonnegative, $w \in L^\infty(0,T),$ and it satisfies
\begin{equation*}
w(t)\leq w_0+a\int^t_0w(s)\ln[a+w(s)]ds, t\in[0,T],
\end{equation*}
where $w_0\geq0,$ $a\geq1.$

Then we have
\begin{equation*}
w(t)\leq (a+w_0)^{e^{at}}, t\in[0,T].
\end{equation*}
\end{lemma}

Next, a difference inequality is given as follows.

\begin{lemma}[\cite{nak1}]\label{dif}
Let $\phi$ be a nonnegative decreasing function on $\mathbb R^+,$ satisfying
\begin{equation*}
  \phi(t+1)-d_2 \phi(t)\leq d_3\big(\phi(t)-\phi(t+1)\big),
\end{equation*}
with some constants $0<d_2<1, d_3>0.$ Then we have
\begin{equation*}
 \displaystyle \phi(t)\leq ce^{-kt},
\end{equation*}
for some constants $c, k>0.$
\end{lemma}

\section{Existence of weak solutions to (\ref{1.1})}

In this section, we prove the existence theorem of weak solutions for system (\ref{1.1}) (Theorem \ref{t1.3}).
\noindent {\bf Proof of Theorem \ref{t1.3}.} For each $\varepsilon>0,$ consider the following perturbation wave equation to equation (\ref{1.1}) in a cylindrical domain $Q_T:$
\begin{equation}\label{1.3}
\left\{\begin{array}{ll}
u^{''}_{\varepsilon}-\Delta u_\varepsilon+au'_\varepsilon+bu_\varepsilon+\dfrac{1}{\varepsilon}\chi u_\varepsilon=u_\varepsilon\ln|u_\varepsilon|^\gamma& (x,t)\in Q_T,\\[2mm]
u_\varepsilon=0 &(x,t)\in \Sigma_T,\\[2mm]
u_\varepsilon(x,0)=\widetilde{u}_0(x),\ u'_\varepsilon(x,0)=\widetilde{u}_1(x) &x\in \Omega,
\end{array}\right.
\end{equation}
where $Q_T=\Omega\times(0,T)\supset\widehat{Q}_T,$ $\Sigma_T$ is its lateral boundary, $\chi$ is a characteristic function defined by
\begin{equation*}
  \chi(x,t)=\left\{
              \begin{array}{ll}
                1 & (x,t)\in Q_T-\widehat{Q}_T, \\
                0 & (x,t)\in\widehat{Q}_T,
              \end{array}
            \right.
\end{equation*}
and the initial data is given by
\begin{equation*}
 \widetilde{ u}_0(x)=\left\{
     \begin{array}{ll}
       u_0(x) & x\in\Omega_0, \\
       0 & x\in\Omega-\Omega_0,
     \end{array}
   \right.
\quad\widetilde{ u}_1(x)=\left\{
     \begin{array}{ll}
       u_1(x) & x\in\Omega_0, \\
       0 & x\in\Omega-\Omega_0.
     \end{array}
   \right.
\end{equation*}

\medskip

We are concerned with  weak solutions of the system (\ref{1.3}).

\medskip

We prove (\ref{1.3}) has a weak solution by Galerkin method.
First, let $\{w_j\}^\infty_{j=1}$ be a orthonormal basis of $H^1_0(\Omega)$ and $V_m$ a subspace spanned by $w_1,w_2,\cdots,w_m.$

\medskip

{\bf 1. An approximate problem} Consider the following approximation problem:

let $u_{\varepsilon m}(t)=\sum\limits^m_{j=1}g_{mj}(t)w_j$ be a solution of ordinary differential equations:
\begin{equation}\label{1.4}
\left\{\begin{array}{ll}
\big(u^{''}_{\varepsilon m}(t),w_j\big)+\big(\nabla u_{\varepsilon m}(t),\nabla w_j\big)+\big(au'_{\varepsilon m}(t)+bu_{\varepsilon m}(t),w_j\big)\\[3mm]
\qquad+\dfrac{1}{\varepsilon }\big(\chi(t) u_{\varepsilon m}(t),w_j\big)=\big(u_{\varepsilon m}(t)\ln|u_{\varepsilon m}(t)|^\gamma,w_j\big), \ j=1,2,\cdots,m,\\[3mm]
u_{\varepsilon m}(0)=\sum\limits^m_{j=1}\big(\widetilde{u}_0,w_j)w_j,\\[3mm]
u'_{\varepsilon m}(0)=\sum\limits^m_{j=1}\big(\widetilde{u}_1,w_j)w_j,
\end{array}\right.
\end{equation}
By the Cauchy-Peano theorem, we know that system (\ref{1.4}) admits a solution $u_{\varepsilon m}$ over $[0,T_m],$ with $g_{mj}\in C^2[0,T_m].$ Next, a priori estimate to be obtained shows that $T_m$ can be replaced by $T.$

\medskip

{\bf 2. A priori estimate}
Multiplying both sides of the first equation in (\ref{1.4})
by $2g'_{mj}$ and adding on $j$ from 1 to $m,$
we obtain
\begin{eqnarray}\label{1.5}
\begin{array}{rl}
&\dfrac{d}{dt}\(\|u'_{\varepsilon m}(t)\|^2+\|\nabla u_{\varepsilon m}(t)\|^2+b\|u_{\varepsilon m}(t)\|^2\)+a\|u'_{\varepsilon m}(t)\|^2\\[3mm]
&\quad +\dfrac{2}{\varepsilon }\big(\chi(t) u_{\varepsilon m}(t),u'_{\varepsilon m}(t)\big)=2\big(u_{\varepsilon m}(t)\ln|u_{\varepsilon m}(t)|^\gamma,u'_{\varepsilon m}(t)\big).
\end{array}
\end{eqnarray}
Set $\displaystyle g_\mu(t)=\frac{1}{\mu}\int^{t+\mu}_{t}\chi(s)ds.$ Then one has
\begin{equation*}
 g_\mu(t)\rightarrow\chi(t),\quad \mbox{with} \quad \mu\rightarrow 0,
\end{equation*}
and by the condition (H1),
\begin{equation*}
 g'_\mu(t)=\dfrac{1}{\mu}\big[\chi(t+\mu)-\chi(t)\big]\leq0.
\end{equation*}
Moreover,
\begin{eqnarray}\label{1.6}
\begin{array}{rl}
\dfrac{2}{\varepsilon }\big(\chi(t) u_{\varepsilon m}(t),u'_{\varepsilon m}(t)\big)\!\!\!&=\lim\limits_{\mu\rightarrow0}\dfrac{2}{\varepsilon }\big(g_\mu(t) u_{\varepsilon m}(t),u'_{\varepsilon m}(t)\big)\\[3mm]
&=\lim\limits_{\mu\rightarrow0}\[\dfrac{d}{dt}\dfrac{1}{\varepsilon}\big(g_\mu(t) u_{\varepsilon m}(t),u_{\varepsilon m}(t)\big)-\dfrac{1}{\varepsilon}\big(g'_\mu(t) u_{\varepsilon m}(t),u_{\varepsilon m}(t)\big)\]\\[3mm]
&\geq\dfrac{d}{dt}\dfrac{1}{\varepsilon }\big(\chi(t) u_{\varepsilon m}(t),u_{\varepsilon m}(t)\big).
\end{array}
\end{eqnarray}
Also, it is easy to check that
\begin{equation}\label{1.7}
2\big(u_{\varepsilon m}(t)\ln|u_{\varepsilon m}(t)|^\gamma,u'_{\varepsilon m}(t)\big)=\dfrac{d}{dt}\big[\big(u_{\varepsilon m}(t)\ln|u_{\varepsilon m}(t)|^\gamma,u_{\varepsilon m}(t)\big)-\dfrac{\gamma}{2}\|u_{\varepsilon m}(t)\|^2\big].
\end{equation}
Define the energy of $u_{\varepsilon m}(t)$ as follows.
\begin{eqnarray*}
\begin{array}{rl}
 E_{\varepsilon m}(t)=E\big(u_{\varepsilon m}(t)\big)=\!\!\!\!&\displaystyle\frac{1}{2}\|u'_{\varepsilon m}(t)\|^2+\frac{1}{2}\|\nabla u_{\varepsilon m}(t)\|^2+\frac{1}{2}\big(b+\frac{\gamma}{2}\big)\|u_{\varepsilon m}(t)\|^2+\frac{1}{2 \varepsilon}\|\chi(t) u_{\varepsilon m}(t)\|^2\\[3mm]
&\displaystyle-\frac{1}{2}\int_\Omega u^2_{\varepsilon m}(t)\ln|u_{\varepsilon m}(t)|^\gamma dx.
\end{array}
\end{eqnarray*}
Substituting (\ref{1.6}) and (\ref{1.7}) into (\ref{1.5}), and integrating it on $(0,t),$ we get
\begin{eqnarray}\label{1.8}
\displaystyle E\big(u_{\varepsilon m}(t)\big)+a\int^t_0\|u'_{\varepsilon m}(\tau)\|^2d\tau\leq E\big(u_{\varepsilon m}(0)\big).
\end{eqnarray}
Since $H^1_0(\Omega)\subset C(\overline{\Omega}),$ $\Omega\subset\mathbb R^1,$ and $u_{\varepsilon m}(0)\rightarrow \widetilde{u}_0$ in $H^1_0(\Omega),$ using the elemental inequality
\begin{equation*}
  \displaystyle|t^2\ln t |\leq C\big(1+t^3\big), \quad \forall t\geq0,
\end{equation*}
we have
\begin{eqnarray}\label{zg00}
\begin{array}{rl}
-\big(u_{\varepsilon m}(0)\ln|u_{\varepsilon m}(0)|^\gamma,u_{\varepsilon m}(0)\big)\!\!\!&\displaystyle\leq\int_\Omega\gamma\big\vert u^2_{\varepsilon m}(0)\ln|u_{\varepsilon m}(0)|\big\vert dx\\[4mm]
&\leq C\big(1+\|\widetilde{u}_0\|^3_{H^1_0(\Omega)}\big).
\end{array}
\end{eqnarray}
In addition, $\forall \varepsilon>0,$
\begin{eqnarray}\label{ff}
\begin{array}{rl}
&\displaystyle\dfrac{1}{\varepsilon }\|\chi(0) u_{\varepsilon m}(0)\|^2\rightarrow0,\quad m\rightarrow\infty.\\[3mm]
\end{array}
\end{eqnarray}
Combining (\ref{zg00}) and (\ref{ff}), we have $\forall \varepsilon>0,$ for large m,
$$\displaystyle E\big(u_{\varepsilon m}(0)\big)\leq C(\widetilde{u}_0, \widetilde{u}_1),$$
where $C$ is a positive constant only depending on the initial data $\widetilde{u}_0, \widetilde{u}_1.$

\medskip

Furthermore, (\ref{1.8}) yields
\begin{equation}\label{e+}
 E\big(u_{\varepsilon m}(t)\big)+a\int^t_0\|u'_{\varepsilon m}(\tau)\|^2d\tau\leq C(\widetilde{u}_0, \widetilde{u}_1).
\end{equation}
Let
\begin{equation*}
E^+\big(u_{\varepsilon m}(t)\big)=\frac{1}{2}\|u'_{\varepsilon m}(t)\|^2+\frac{1}{2}\|\nabla u_{\varepsilon m}(t)\|^2+\frac{1}{2}\big(b+\dfrac{\gamma}{2}\big)\|u_{\varepsilon m}(t)\|^2+\dfrac{1}{2\varepsilon }\|\chi(t) u_{\varepsilon m}(t)\|^2.
\end{equation*}
Then (\ref{e+}) becomes
\begin{equation}\label{e}
 E^+\big(u_{\varepsilon m}(t)\big)-\frac{1}{2}\big(u_{\varepsilon m}(t)\ln|u_{\varepsilon m}(t)|^\gamma,u_{\varepsilon m}(t)\big)+a\int^t_0\|u'_{\varepsilon m}(\tau)\|^2d\tau\leq C(\widetilde{u}_0, \widetilde{u}_1).
\end{equation}
Using the logarithmic Sobolev inequality in (\ref{e}), we get $\forall\tilde{a}\in\mathbb R,$
\begin{eqnarray*}\label{}
\begin{array}{rl}
&\displaystyle E^+\big(u_{\varepsilon m}(t)\big)-\frac{\gamma\widetilde{a}^2}{4\pi}\|\nabla u_{\varepsilon m}(t)\|^2+\frac{\gamma (1+\ln\tilde{a})}{4}\|u_{\varepsilon m}(t)\|^2+a\int^t_0\|u'_{\varepsilon m}(\tau)\|^2d\tau\\[3mm]
&\displaystyle\leq C(\widetilde{u}_0, \widetilde{u}_1)+\frac{\gamma}{2}\|u_{\varepsilon m}(t)\|^2\ln\| u_{\varepsilon m}(t)\|.
\end{array}
\end{eqnarray*}
Precisely,
\begin{eqnarray*}\label{}
\begin{array}{rl}
&\displaystyle \frac{1}{2}\|u'_{\varepsilon m}(t)\|^2+(\frac{1}{2}-\frac{\gamma\widetilde{a}^2}{4\pi})\|\nabla u_{\varepsilon m}(t)\|^2+\big[\frac{1}{2}(b+\frac{\gamma}{2})+\frac{\gamma (1+\ln\tilde{a})}{4}\big]\|u_{\varepsilon m}(t)\|^2\\[3mm]
&\displaystyle+a\int^t_0\|u'_{\varepsilon m}(\tau)\|^2d\tau\leq C(\widetilde{u}_0, \widetilde{u}_1)+\frac{\gamma}{2}\|u_{\varepsilon m}(t)\|^2\ln\| u_{\varepsilon m}(t)\|.
\end{array}
\end{eqnarray*}
Without loss of generality, let $\displaystyle\widetilde{a}=(\frac{\pi}{\gamma})^\frac{1}{2},$ then $\displaystyle\frac{1}{2}-\frac{\gamma\widetilde{a}^2}{4\pi}=\frac{1}{4}$ and $\displaystyle\big[\frac{1}{2}(b+\frac{\gamma}{2})+\frac{\gamma (1+\ln\tilde{a})}{4}\big]>0$ ($0<\gamma<1$). Consequently, we arrive at
\begin{eqnarray}\label{11}
\begin{array}{rl}
&\displaystyle\|u'_{\varepsilon m}(t)\|^2+\|\nabla u_{\varepsilon m}(t)\|^2+\|u_{\varepsilon m}(t)\|^2+\dfrac{1}{\varepsilon }\|\chi(t) u_{\varepsilon m}(t)\|^2+2a\int^t_0\|u'_{\varepsilon m}(\tau)\|^2d\tau\\[3mm]
&\leq C_1+C_2\| u_{\varepsilon m}(t)\|^2\ln\| u_{\varepsilon m}(t)\|^2,
\end{array}
\end{eqnarray}
where $C_i(i=1,2)$ are positive constants just related to initial data $\widetilde{u}_0, \widetilde{u}_1$ and known parameters.

\medskip

As
\begin{equation*}
 \displaystyle u_{\varepsilon m}(t)=u_{\varepsilon m}(0)+\int^t_0\frac{\partial u_{\varepsilon m}}{\partial\tau}(\tau)d\tau,
\end{equation*}
one has
\begin{eqnarray*}
&&|u_{\varepsilon m}(t)|^2\leq2|u_{\varepsilon m}(0)|^2+
2t\int^t_0|\frac{\partial u_{\varepsilon m}}{\partial\tau}(\tau)|^2d\tau.
\end{eqnarray*}
Integrating above inequality on $\Omega,$ we have
\begin{equation}\label{zg3}
\|u_{\varepsilon m}(t)\|^2\leq2\|\widetilde{u}_0\|^2+\max\{2T,1\}\big(\frac{1+C_2}{C_2}\big)
\int^t_0\|\frac{\partial u_{\varepsilon m}}{\partial\tau}(\tau)\|^2d\tau.
\end{equation}
 (\ref{11}) together with (\ref{zg3}), yield
\begin{eqnarray}\label{zg1}
\begin{array}{rl}
\|u_{\varepsilon m}(t)\|^2\!\!\!&\displaystyle\leq2\|\widetilde{u}_0\|^2+\max\{2T,1\}\big(\frac{1+C_2}{C_2}\big)C_1T\\[3mm]
&\displaystyle\quad+\max\{2T,1\}\big(1+C_2\big)
\int^t_0\| u_{\varepsilon m}(\tau)\|^2\ln\| u_{\varepsilon m}(\tau)\|^2d\tau.
\end{array}
\end{eqnarray}
Put $\displaystyle A=2\|\widetilde{u}_0\|^2+\max\{2T,1\}\big(\frac{1+C_2}{C_2}\big)C_1T$ and $\displaystyle B=\max\{2T,1\}\big(1+C_2\big).$ According to the logarithmic Gronwall inequality, we derive
\begin{eqnarray}\label{1.11}
\|u_{\varepsilon m}(t)\|^2\leq\big(A+B\big)^{e^{Bt}}\leq C_T.
\end{eqnarray}
Applying (\ref{1.11}) into (\ref{11}), we obtain that for $a.e.\ t\in[0,T],$
\begin{eqnarray}\label{1.12}
&\displaystyle\|u'_{\varepsilon m}(t)\|^2+\|\nabla u_{\varepsilon m}(t)\|^2+\|u_{\varepsilon m}(t)\|^2+\dfrac{1}{\varepsilon }\|\chi(t) u_{\varepsilon m}(t)\|^2+2a\int^t_0\|u'_{\varepsilon m}(\tau)\|^2d\tau\leq C_T.
\end{eqnarray}
The estimate (\ref{1.12}) implies that $T_k=T.$

\medskip

{\bf 3. Passing to the limit }
Until now, we know that the sequence of approximate solutions $\big(u_{\varepsilon m}\big)_m
$ to (\ref{1.4}) satisfy
\begin{equation}\label{1.13}
\begin{array}{rl}
&\displaystyle\big(u_{\varepsilon m}\big)_m \ \mbox{is a bounded sequence in}\ L^\infty\big(0,T; H^1_0(\Omega)\big),\\[3mm]
&\displaystyle\big(u'_{\varepsilon m}\big)_m \ \mbox{is a bounded sequence in}\ L^\infty\big(0,T; L^2(\Omega)\big),\\[3mm]
&\displaystyle\big(\frac{1}{\sqrt{\varepsilon}}\chi u_{\varepsilon m}\big)_m \ \mbox{is a bounded sequence in}\ L^\infty\big(0,T; L^2(\Omega)\big).
\end{array}
\end{equation}
We extract a subsequence, still written as $\big(u_{\varepsilon m}\big)_m,$ and there exists a $u_{\varepsilon }$ such that
\begin{equation}\label{1.14}
\begin{array}{rl}
&\displaystyle u_{\varepsilon m}\stackrel{*}\rightharpoonup u_{\varepsilon }\ \mbox{in}\ L^\infty\big(0,T; H^1_0(\Omega)\big),\\[3mm]
&\displaystyle  u'_{\varepsilon m}\stackrel{*}\rightharpoonup u'_{\varepsilon }\ \mbox{ in}\ L^\infty\big(0,T; L^2(\Omega)\big),\\[3mm]
&\displaystyle u_{\varepsilon m}\rightharpoonup u_{\varepsilon }\ \mbox{in}\ L^2\big(0,T; H^1_0(\Omega)\big),\\[3mm]
&\displaystyle  u'_{\varepsilon m}\rightharpoonup u'_{\varepsilon }\ \mbox{ in}\ L^2\big(0,T; L^2(\Omega)\big).
\end{array}
\end{equation}
By Aubin-Lions lemma, we also obtain the following convergence from (\ref{1.13}).
\begin{equation}\label{1.15}
\displaystyle u_{\varepsilon m}\rightarrow u_{\varepsilon }\ \mbox{in}\ L^2\big(0,T; L^2(\Omega)\big).
\end{equation}
So there exists a subsequence
\begin{equation*}
u_{\varepsilon m }\rightarrow u_{\varepsilon }\quad \ a.e. \ \mbox{in} \ [0,T]\times \Omega.
\end{equation*}
Since the mapping $x\mapsto x\ln|x|^\gamma$ is continuous,
\begin{equation*}
u_{\varepsilon m }(t)\ln|u_{\varepsilon m}(t)|^\gamma\rightarrow u_{\varepsilon }(t)\ln|u_{\varepsilon }(t)|^\gamma\quad \ a.e. \ \mbox{in} \ [0,T]\times \Omega.
\end{equation*}
By Lebesgue dominated convergence theorem, we obtain
\begin{equation}\label{1.16}
u_{\varepsilon m }(t)\ln|u_{\varepsilon m}(t)|^\gamma\rightarrow u_{\varepsilon }(t)\ln|u_{\varepsilon }(t)|^\gamma\quad \mbox{in} \ L^2\big(0,T;L^2(\Omega)\big).
\end{equation}
Take (\ref{1.14}), (\ref{1.15}) and (\ref{1.16}) into (\ref{1.4}), then we have for any $\theta\in C^\infty_0(0,T),$
\begin{equation}\label{1.17}
\left\{\begin{array}{ll}
\displaystyle-\int^T_0\big(u^{'}_{\varepsilon}(t),w_j \theta'(t)\big)dt+\int^T_0\big(\nabla u_{\varepsilon }(t),\nabla w_j\theta(t)\big)dt+\int^T_0\big(au'_{\varepsilon }(t)+bu_{\varepsilon }(t),w_j\theta(t)\big)dt\\[3mm]
\displaystyle\qquad+\int^T_0\dfrac{1}{\varepsilon }\big(\chi(t) u_{\varepsilon }(t),w_j\theta(t)\big)dt=\int^T_0\big(u_{\varepsilon }(t)\ln|u_{\varepsilon }(t)|^\gamma,w_j\theta(t)\big)dt, \\[3mm]
\displaystyle u_{\varepsilon }(0)=\widetilde{u}_0,\\[3mm]
u'_{\varepsilon }(0)=\widetilde{u}_1,
\end{array}\right.
\end{equation}
which  means that $u_{\varepsilon }$ is a weak solution of (\ref{1.3}).

\medskip

One can check  the bounds (\ref{1.13}) obtained on $u_{\varepsilon m}$  are also valid for $u_{\varepsilon }.$ With a similar argument, we get
\begin{equation*}\label{1.18}
\begin{array}{rl}
&\displaystyle u_{\varepsilon }\stackrel{*}\rightharpoonup u\ \mbox{in}\ L^\infty\big(0,T; H^1_0(\Omega)\big),\\[3mm]
&\displaystyle  u'_{\varepsilon }\stackrel{*}\rightharpoonup u'\ \mbox{in}\ L^\infty\big(0,T; L^2(\Omega)\big),\\[3mm]
&\displaystyle u_{\varepsilon }\rightharpoonup u\ \mbox{in}\ L^2\big(0,T; H^1_0(\Omega)\big),\\[3mm]
&\displaystyle  u'_{\varepsilon }\rightharpoonup u'\ \mbox{in}\ L^2\big(0,T; L^2(\Omega)\big).
\end{array}
\end{equation*}
Aubin-Lions lemma also tell us
\begin{equation}\label{1.19}
\displaystyle u_{\varepsilon }\rightarrow u\ \mbox{in}\ L^2\big(0,T; L^2(\Omega)\big).
\end{equation}
As $\displaystyle\frac{1}{\sqrt{\varepsilon}}\chi u_{\varepsilon }$ remains bounded in $L^\infty\big(0,T; L^2(\Omega)\big),$ we have $\|\chi u_{\varepsilon }\|^2\leq C\varepsilon$ and affirm that
\begin{equation*}
\displaystyle\chi u_{\varepsilon }\rightarrow 0 \quad \mbox{in} \ L^2\big(0,T; L^2(\Omega)\big).
\end{equation*}
It follows from (\ref{1.19}) that
\begin{equation*}
\chi u=0 \quad a.e. \ \mbox{in} \ [0,T]\times \Omega.
\end{equation*}
Since $\chi=1$ on $Q_T-\widehat{Q}_T,$ we deduce
\begin{equation*}\label{1.20}
u=0\quad a.e. \ \mbox{in} \ Q_T-\widehat{Q}_T.
\end{equation*}
Observe that linear combinations of functions of the form $\theta\otimes w_j$ are dense in the space $C^\infty_0\big(0,T; H^1_0(\Omega)\big).$ Let $v$ be the restriction of $u$ on $\widehat{Q}_T.$ Then, for any $\phi\in C_0^\infty\big(0,T; H^1_0(\Omega_t)\big)\subset C_0^\infty\big(0,T; H^1_0(\Omega)\big),$ it holds that
\begin{equation}\label{1.21}
\begin{array}{ll}
&\displaystyle-\int_{\widehat{Q}_T}v^{'}(t)\phi'(t)dxdt+\int_{\widehat{Q}_T}\nabla v(t)\nabla \phi(t)dxdt+\int_{\widehat{Q}_T}\big(av'(t)+bv(t)\big)\phi(t)dxdt\\[3mm]
&=\displaystyle\int_{\widehat{Q}_T}v(t)\ln|v(t)|^\gamma\phi(t)dxdt,
\end{array}
\end{equation}
(\ref{1.21}) indicates $v$ is a weak solution of system (\ref{1.1}).

\medskip

From the fact $u(t)\in H^1_0(\Omega)$ and the assumption (H2), we get
\begin{eqnarray*}\label{1.22}
v\in L^\infty\big(0,T; H^1_0(\Omega_t)\big),
\end{eqnarray*}
Moreover,
\begin{eqnarray*}\label{1.23}
v'\in L^\infty\big(0,T; L^2(\Omega_t)\big).
\end{eqnarray*}
Recall that supp$\widetilde{u}_0\subset\Omega_0.$ Hence, $v$ satisfies the initial condition
\begin{eqnarray*}\label{1.24}
v(0)=u_0, \quad v'(0)=u_1.
\end{eqnarray*}
In fact, $v$ satisfies
\begin{equation*}\label{1.25}
\begin{array}{ll}
&\displaystyle-\int_{\widehat{Q}_T}v^{'}(t)\varphi'(t)dxdt+\int_{\widehat{Q}_T}\nabla v(t)\cdot\nabla \varphi(t)dxdt+\int_{\widehat{Q}_T}\big(av'(t)+bv(t)\big)\varphi(t)dxdt\\[3mm]
&\displaystyle=\int_{\widehat{Q}_T}v(t)\ln|v(t)|^\gamma\varphi(t)dxdt+\big(u_1,\varphi(0)\big)_{\Omega_0},
\end{array}
\end{equation*}
for all $\displaystyle \varphi\in L^1\big(0,T; H^1_0(\Omega_t)\big),$ $\displaystyle \varphi'\in L^1\big(0,T; L^2(\Omega_t)\big)$ with $\varphi(x,T)=0.$
\endpf

\section{Decay of solutions to (\ref{1.1})}

In this section, we study the decay property of solutions to (\ref{1.1}).

\medskip

Let
\begin{equation*}
\begin{array}{rl}
 J_{\varepsilon m}(t)=J\big(u_{\varepsilon m}(t)\big)=\!\!\!&\displaystyle\frac{1}{2}\|\nabla u_{\varepsilon m}(t)\|^2+\frac{1}{2}\big(b+\frac{\gamma}{2}\big)\|u_{\varepsilon m}(t)\|^2+\frac{1}{2 \varepsilon}\|\chi(t) u_{\varepsilon m}(t)\|^2\\[3mm]
&\displaystyle-\frac{1}{2}\int_\Omega u^2_{\varepsilon m}(t)\ln|u_{\varepsilon m}(t)|^\gamma dx.
\end{array}
\end{equation*}
Hence,
\begin{equation*}
  E_{\varepsilon m}(t)=\frac{1}{2}\|u'_{\varepsilon m}(t)\|^2+ J_{\varepsilon m}(t).
\end{equation*}
We call $J_{\varepsilon m}(t)$ the potential of $u_{\varepsilon m}$ at time $t.$

\medskip

Firstly, we give the the following lemma which guarantee the energy is positive.
\begin{lemma}\label{pos}
If $\widetilde{u}_0\in H^1_0(\Omega),$ and $\widetilde{u}_1\in L^2(\Omega),$ such that $0<E(0)<d,$ and $I_1(\widetilde{u}_0)>0,$ then $\forall \varepsilon>0,$
\begin{equation*}
u_{\varepsilon m}(t)\in W:=\{J\big(u_{\varepsilon m}(t)\big)<d,\ I_1\big(u_{\varepsilon m}(t)\big)>0\}\cup\{0\},
\end{equation*}
for $t\in[0,T]$ and large $m.$
\end{lemma}

\noindent {\bf Proof. }
It is easy to check that
\begin{equation*}
 \displaystyle \frac{d E_{\varepsilon m}(t)}{dt}\leq-a\|u'_{\varepsilon m}(t)\|^2<0.
\end{equation*}
This implies that for any fixed $\varepsilon>0,$ there exists a integer $m_0,$ such that for $t\in[0,T],$
\begin{equation}\label{d1.1}
J\big(u_{\varepsilon m}(t)\big)+\frac{1}{2}\|u'_{\varepsilon m}(t)\|^2=E_{\varepsilon m}(t)\leq E_{\varepsilon m}(0)<d, \quad m\geq m_0.
\end{equation}
Thus, in this case,
\begin{equation*}
 J\big(u_{\varepsilon m}(t)\big)<d.
\end{equation*}
As
\begin{equation*}
I_1\big(\widetilde{u}_0\big)=\|\nabla \widetilde{u}_0\|^2-\int_\Omega \widetilde{u}_0^2\ln|\widetilde{u}_0|^\gamma dx>0,
\end{equation*}
and $\displaystyle u_{\varepsilon m }(0)\rightarrow \widetilde{u}_0$ in $H^1_0(\Omega),$ we have $I_1\big(u_{\varepsilon m }(0)\big)>0,$ $m\geq m_0.$ Next, we claim that
\begin{equation*}
u_{\varepsilon m}(t)\in W,\quad \forall t\in[0,T],\ m\geq m_0.
\end{equation*}
Otherwise, there exist a point $t_1\in[0,T]$ and a integer $m_1\geq m_0,$ such that $u_{\varepsilon m }(t_1)\in\partial W.$ Thus $(1)$ $J\big(u_{\varepsilon m_1}(t_1)\big)=d$ (impossible).
\begin{equation*}
(2)\quad I_1\big(u_{\varepsilon m_1 }(t_1)\big)=0,\quad \mbox{and} \quad \|\nabla u_{\varepsilon m_1 }(t_1)\|\neq0.
\end{equation*}
Then $u_{\varepsilon m_1 }(t_1)\in N,$ and
\begin{equation*}
 J\big(u_{\varepsilon m_1 }(t_1)\big)\geq J_1\big(u_{\varepsilon m_1 }(t_1)\big)\geq\inf\limits_{u\in N}J_1(u)=d.
\end{equation*}
However, according to (\ref{d1.1}), we know that $J\big(u_{\varepsilon m_1 }(t_1)\big)<d,$ a contradiction.
\endpf

\medskip

\noindent {\bf Proof of Theorem \ref{t1}.} For simplicity, let $a=1$ and $b=1.$ To finish the proof of Theorem \ref{t1.3}, it suffices to show that the approximate solutions  $u_{\varepsilon m}$ (m: large) satisfy decay estimate independent of $ \varepsilon, m.$ Replacing $w_j$ by $u'_{\varepsilon m}$ in the first equation of (\ref{1.4}) and integrating it on $(t_1,t_2)\times \Omega,$ we obtain
\begin{eqnarray}\label{p1}
\begin{array}{rl}
\displaystyle E_{\varepsilon m}(t_2)+\int^{t_2}_{t_1}\|u'_{\varepsilon m}(\tau)\|^2d\tau\leq E_{\varepsilon m}(t_1).
\end{array}
\end{eqnarray}
Notice that $E_{\varepsilon m}$ is a nonnegative monotone decreasing function on $\mathbb R^+.$ By (\ref{p1}) with $t_1=t,$ $t_2=t+1,$ we see that there exists two points $\displaystyle t_1\in [t,t+\frac{1}{4}]$ and $\displaystyle t_2\in [t+\frac{3}{4}, t+1]$ such that
\begin{eqnarray}\label{p2}
\begin{array}{rl}
\|u'_{\varepsilon m}(t_i)\|\leq2D_{\varepsilon m}(t),
\end{array}
\end{eqnarray}
where $D^2_{\varepsilon m}(t)=E_{\varepsilon m}(t)-E_{\varepsilon m}(t+1).$

\medskip

Replacing $w_j$ by $u_{\varepsilon m}$ in the first equation of (\ref{1.4}) and integrating it on $(t_1,t_2)\times \Omega,$ we obtain
\begin{eqnarray}\label{p3}
&&\displaystyle\int_{t_1}^{t_2}\int_\Omega \Big(
|\nabla u_{\varepsilon m}|^2+u^2_{\varepsilon m}+\frac{1}{\varepsilon}\chi u^2_{\varepsilon m}-u^2_{\varepsilon m}\ln|u_{\varepsilon m}|^\gamma\Big)dxd\tau
+\displaystyle\int_{t_1}^{t_2}\int_\Omega u'_{\varepsilon m}u_{\varepsilon m}dxd\tau\nonumber\\[3mm]
&&= \big( u'_{\varepsilon m}(t_1),u_{\varepsilon m}(t_1)\big)
-\big( u'_{\varepsilon m}(t_2),u_{\varepsilon m}(t_2)\big)
+\displaystyle\int_{t_1}^{t_2}\| u'_{\varepsilon m}(\tau)\|^2d\tau.
\end{eqnarray}
Using Cauchy inequality, from (\ref{p3}) we deduce
\begin{eqnarray*}\label{}
&&\displaystyle\int_{t_1}^{t_2}\int_\Omega \Big(
|\nabla u_{\varepsilon m}|^2+\frac{1}{2}\big(1+\frac{\gamma}{2}\big)u^2_{\varepsilon m}+\frac{1}{\varepsilon}\chi u^2_{\varepsilon m}-u^2_{\varepsilon m}\ln|u_{\varepsilon m}|^\gamma\Big)dxd\tau\nonumber\\[3mm]
&&\leq \big( u'_{\varepsilon m}(t_1),u_{\varepsilon m}(t_1)\big)
-\big( u'_{\varepsilon m}(t_2),u_{\varepsilon m}(t_2)\big)
+\displaystyle\big(1+\frac{1}{2-\gamma}\big)\int_{t_1}^{t_2}\| u'_{\varepsilon m}(\tau)\|^2d\tau.
\end{eqnarray*}
It is known from Lemma \ref{pos} that $\displaystyle I_1\big(u_{\varepsilon m}(t)\big)=\|\nabla u_{\varepsilon m}(t)\|^2-\int_\Omega u_{\varepsilon m}^2(t)\ln|u_{\varepsilon m}(t)|^\gamma dx>0.$ By (\ref{p2}) we get further
{\setlength\arraycolsep{2pt}
\begin{eqnarray}\label{p4}
&&\displaystyle\int_{t_1}^{t_2}\int_\Omega \Big(\frac{1}{2}
|\nabla u_{\varepsilon m}|^2+\frac{1}{2}\big(1+\frac{\gamma}{2}\big)u^2_{\varepsilon m}+\frac{1}{2\varepsilon}\chi u^2_{\varepsilon m}-\frac{1}{2}u^2_{\varepsilon m}\ln|u_{\varepsilon m}|^\gamma\Big)dxd\tau\nonumber\\[3mm]
&&\leq 2D_{\varepsilon m}(t)\Big(\|u_{\varepsilon m}(t_1)\|+\|u_{\varepsilon m}(t_2)\|\Big)
+\displaystyle\big(1+\frac{1}{2-\gamma}\big)\int_{t_1}^{t_2}\| u'_{\varepsilon m}(\tau)\|^2d\tau.
\end{eqnarray}}
Combining (\ref{p1}) and (\ref{p4}), we arrive at
\begin{eqnarray*}
\begin{array}{rl}
 \displaystyle\int^{t_2}_{t_1}E_{\varepsilon m}(\tau)d\tau\!\!\!&\leq2D_{\varepsilon m}(t)\Big(\|u_{\varepsilon m}(t_1)\|+\|u_{\varepsilon m}(t_2)\|\Big)
+\displaystyle\big(\frac{3}{2}+\frac{1}{2-\gamma}\big)\int_{t_1}^{t_2}\| u'_{\varepsilon m}(\tau)\|^2d\tau\\[4mm]
&\displaystyle\leq2D_{\varepsilon m}(t)\Big(\|u_{\varepsilon m}(t_1)\|+\|u_{\varepsilon m}(t_2)\|\Big)+\big(\frac{3}{2}+\frac{1}{2-\gamma}\big)D^2_{\varepsilon m}(t).
\end{array}
\end{eqnarray*}
Since $\displaystyle\|u_{\varepsilon m}(t_i)\|^2\leq2E_{\varepsilon m}(t_i)\leq2E_{\varepsilon m}(t),$ $(i=1,2),$ using Cauchy inequality again, we obtain for any $\delta>0,$
\begin{eqnarray*}
\displaystyle \int^{t_2}_{t_1}E_{\varepsilon m}(\tau)d\tau\!\!\!\!\!\!\!\!&&\displaystyle \leq\frac{1}{\delta}D^2_{\varepsilon m}(t)+\delta\Big(\|u_{\varepsilon m}(t_1)\|^2+\|u_{\varepsilon m}(t_2)\|^2\Big)
+\displaystyle\big(\frac{3}{2}+\frac{1}{2-\gamma}\big)D^2_{\varepsilon m}(t)\\[3mm]
&& \displaystyle \leq\frac{1}{\delta}D^2_{\varepsilon m}(t)+4\delta E_{\varepsilon m}(t)+\big(\frac{3}{2}+\frac{1}{2-\gamma}\big)D^2_{\varepsilon m}(t).
\end{eqnarray*}
On the other hand, there exists a point $\displaystyle t^*\in[t_1,t_2]\subset[t,t+1],$ such that
\begin{equation*}
\frac{1}{t_2-t_1}\int^{t_2}_{t_1}E_{\varepsilon m}(\tau)d\tau=E_{\varepsilon m}(t^*).
\end{equation*}
Notice that $\displaystyle\frac{1}{2}\leq|t_2-t_1|\leq1,$ and $E_{\varepsilon m}$ is decreasing, so we have
\begin{equation*}
E_{\varepsilon m}(t+1)\leq E_{\varepsilon m}(t^*)\leq\frac{2}{\delta}D^2_{\varepsilon m}(t)+8\delta E_{\varepsilon m}(t)+\big(3+\frac{2}{2-\gamma}\big)D^2_{\varepsilon m}(t).
\end{equation*}
Then
\begin{equation*}
E_{\varepsilon m}(t+1)-8\delta E_{\varepsilon m}(t)\leq\big(\frac{2}{\delta}+3+\frac{2}{2-\gamma}\big)D^2_{\varepsilon m}(t)=\big(\frac{2}{\delta}+3+\frac{2}{2-\gamma}\big)\big(E_{\varepsilon m}(t)-E_{\varepsilon m}(t+1)\big).
\end{equation*}
Let $\displaystyle d_2=8\delta<1,$ $\displaystyle d_3=\frac{2}{\delta}+3+\frac{2}{2-\gamma}.$ By Lemma \ref{dif}, we obtain for some $C,\beta>0,$
\begin{equation}\label{2.5}
  E_{\varepsilon m}(t)\leq Ce^{-\beta t},\quad t>0.
\end{equation}
As we see (\ref{2.5}) holds independent of $\varepsilon,m,$ by Banach-Steinhaus theorem we can conclude that
$$E(t)\leq Ce^{-\beta t},\quad t>0.$$
\endpf

\begin{remark}\label{l1.2}
In this remark, we give a description on positive constants $C,\beta$ in (\ref{2.5}). Let $t$ be arbitrary. If the following inequality holds.
$$E_{\varepsilon m}(t+1)-d_2E_{\varepsilon m}(t)\leq d_3\big(E_{\varepsilon m}(t)-E_{\varepsilon m}(t+1)\big),$$
where $\displaystyle0<d_2=8\delta<1,$ $\displaystyle d_3=\frac{2}{\delta}+3+\frac{2}{2-\gamma}>0.$

\medskip

Then we have
\begin{equation*}
\begin{array}{ll}
\displaystyle E_{\varepsilon m}(t+1)\leq\Big(\frac{d_3+d_2}{d_3+1}\Big)E_{\varepsilon m}(t).
\end{array}
\end{equation*}
where $\displaystyle\frac{d_3+d_2}{d_3+1}<1.$

\medskip

This means that
\begin{equation*}
 \displaystyle E_{\varepsilon m}(t)\leq E_{\varepsilon m}(0)e^{\displaystyle\ln\big(\displaystyle\frac{d_3+d_2}{d_3+1}\big)t},
\end{equation*}
and $\displaystyle\frac{d_3+d_2}{d_3+1}$ gets the minimum at  $\delta=\displaystyle\frac{16-\sqrt{(16)^2+4A}}{-2A},$ with $\displaystyle A=4\big(4+\frac{2}{2-\gamma}\big).$
\end{remark}

\end{document}